\documentclass[11pt]{article}
\usepackage{a4wide}
\usepackage{amsmath,amssymb,amsthm, graphicx}
\usepackage[dvips]{epsfig}
\usepackage{amsfonts}
\usepackage{enumerate}
\usepackage[usenames,dvipsnames]{xcolor}
\usepackage{bm} 
\usepackage[a4paper,left=25mm]{geometry}
\newtheorem{theorem}{Theorem}[section]
\newtheorem{proposition}[theorem]{Proposition}

\newtheorem{claim}[theorem]{Claim} 
\theoremstyle{remark}

\theoremstyle{definition}

\newcommand{\FF}{\mathbb{F}}
\newcommand{\R}{\mathbb{R}}

\newcommand{\N}{{\mathbb N}}

\newcommand{\F}{\mathcal{F}}

\newif\ifcomments
\commentstrue     

\ifcomments
\newcommand{\blue}[1]{{\color{blue}\bf #1}}
\newcommand{\red}[1]{{\color{red}\bf #1}}
\newcommand{\gray}[1]{{\color{gray}\bf #1}}

\newcounter{sideremark}
\newcommand{\marrow}{\stepcounter{sideremark}\marginpar{$\boldsymbol{\longleftarrow\scriptstyle\arabic{sideremark}}$}}
\newcommand{\martin}[1]{{\vskip 5pt\textsf{\blue{*** (Martin) \marrow #1}\vskip
5pt}}}
\newcommand{\eran}[1]{{\vskip 5pt\textsf{\blue{*** (Eran) \marrow #1}\vskip
5pt}}}
\else
\newcommand{\blue}[1]{#1}
\newcommand{\red}[1]{}
\newcommand{\gray}[1]{}
\newcommand{\martin}[1]{}
\newcommand{\eran}[1]{}
\fi

\DeclareMathOperator{\conv}{conv}

\DeclareMathOperator{\relint}{relint}

\DeclareMathOperator{\aff}{aff}

\newcommand{\thh}{\tilde{H}}

\newcommand{\ttop}{\tau}
\newcommand{\rk}{\text{rk}}
\newcommand{\zht}{\widehat{0}}
\newcommand{\oht}{\widehat{1}}

\newcommand{\til}{\tilde{L}}
\newcommand{\enp}{\begin{flushright} $\Box$ \end{flushright}}

\begin{document}

\title{Pach's selection theorem does not admit a topological extension}


\author{Imre B\'ar\'any\thanks{R\'{e}nyi Institute, Hungarian Academy of Sciences,
POB 127, 1364 Budapest, Hungary and Department of Mathematics, University College London, Gower Street, London, WC1E 6BT, UK. email: barany@renyi.hu} \and Roy Meshulam\thanks{Department of Mathematics,  Technion - Israel Institute of Technology, Haifa 32000, Israel. 
email: meshulam@math.technion.ac.il} \and Eran Nevo\thanks{Einstein Institute of Mathematics, The Hebrew University of Jerusalem, Jerusalem 91904, Israel. email: nevo@math.huji.ac.il} \and 
Martin Tancer\thanks{Department of Applied Mathematics, Charles University in Prague,
Malostransk\'{e} n\'{a}m\v{e}st\'{\i} 25, 118 00, Praha 1, Czech Republic. email: tancer@kam.mff.cuni.cz}}

\maketitle

\begin{abstract}
Let $U_1,\ldots,U_{d+1}$ be $n$-element sets in $\R^d$.
Pach's selection theorem says that there exist subsets $Z_1\subset U_1,\ldots,Z_{d+1}\subset U_{d+1}$ and
a point $u \in \R^d$ such that each $|Z_i|\ge c_1(d)n$ and $u \in \conv\{
z_1,\ldots,z_{d+1}\}$ for every choice of $z_1 \in Z_1,\ldots,z_{d+1} \in
Z_{d+1}$.
Here we show that this theorem does not admit a topological extension
with linear size sets $Z_i$. However, there is a topological extension
where each $|Z_i|$ is of order $(\log n)^{1/d}$.
\end{abstract}

\section{Introduction}
\label{s:intro}
Pach's homogeneous selection theorem is the following key result in discrete geometry.
\begin{theorem}[Pach \cite{Pach98}]
\label{t:pach}
For $d \geq 1$ there exists a constant $c_1(d)>0$ such that the following holds.
For any $n$-element sets $U_1,\ldots,U_{d+1}$ in $\R^d$, there exist subsets $Z_1\subset U_1,\ldots,Z_{d+1}\subset U_{d+1}$ and a point $u \in \R^d$ such that each $|Z_i|\ge c_1(d)n$ and $u \in \conv \{z_1,\ldots,z_{d+1}\}$ for every choice of $z_1 \in Z_1,\ldots,z_{d+1} \in Z_{d+1}$.
\end{theorem}

This result was proved by B\'{a}r\'{a}ny, F\"{u}redi, and Lov\'{a}sz
\cite{BFL90} for $d=2$ and by Pach \cite{Pach98} for general $d$. Here we show
that this theorem does not admit a topological extension when the size of the
$Z_i$ is linear in $n$, but does admit one when the sizes are of order $(\log
n)^{1/d}$. Now we reformulate
Theorem~\ref{t:pach} and then we state the topological extension.

Throughout the paper we will identify an abstract simplicial complex $X$ with
its geometric realization. For $k \geq 0$, let $X^{(k)}$ denote the $k$-dimensional skeleton of $X$ and let $X(k)$ be the family of $k$-dimensional faces of $X$.
For an abstract simplex
$\sigma=\{v_0,\ldots,v_k\} \in X(k)$, we write $\langle v_0,\ldots,v_k \rangle$
for its geometric realization.

Let $\Delta_{n-1}$ denote the $(n-1)$-simplex. Consider $d+1$ sets $V_1,\ldots,V_{d+1}$, each of size $n$, and their
join
$$(\Delta_{n-1}^{(0)})^{*(d+1)} \cong V_1*\cdots*V_{d+1} := \{\sigma \subset
\bigcup_{i=1}^{d+1} V_i: |\sigma \cap V_i| \leq 1 \text{~for~all~}1 \leq i \leq
d+1 \}.$$

Trivially, there is an affine map $f:(\Delta_{n-1}^{(0)})^{*(d+1)} \to \R^d$ that is a
bijection between $V_i$ and $U_i$ for each $i$ (where $U_i$ are the sets
from the statement of Pach's theorem). In this setting the homogeneous
selection theorem says that there exist subsets $Z_i \subset V_i$ such that
$|Z_i|\ge c_1(d)n$ and

\begin{equation*}
  \label{e:deft}
  \bigcap_{z_1\in Z_1,\ldots, z_{d+1} \in Z_{d+1}} f(\langle z_1,\ldots,z_{d+1}
  \rangle) \neq \emptyset.
\end{equation*}

Assume now that $f$ is not affine but only continuous. For a mapping
$f:(\Delta_{n-1}^{(0)})^{*(d+1)} \rightarrow \R^d$, let $\tau(f)$ denote the
maximal $m$ such that there exist $m$-element subsets $Z_1 \subset V_1,\ldots,
Z_{d+1} \subset V_{d+1}$ that satisfy
\begin{equation*}
  \label{e:deft}
  \bigcap_{z_1\in Z_1,\ldots, z_{d+1} \in Z_{d+1}} f(\langle
  z_1,\ldots,z_{d+1}\rangle) \neq \emptyset.
\end{equation*}

%
%
%

Define the {\it topological Pach number} $\ttop(d,n)$ to be the minimum of $\tau(f)$ as $f$ ranges over all continuous maps
from $(\Delta_{n-1}^{(0)})^{*(d+1)}$ to $\R^d$.
Our main result is the following:
\begin{theorem}
\label{t:ntp}
For $d \geq 1$ there exists a constant $c_2(d)=O(d)$ such that $\ttop(d,n) \leq c_2(d)n^{1/d}$ for all $n \geq (2d)^d$.
\end{theorem}

For a lower bound on $\ttop(d,n)$ we only have the following:
\begin{theorem}
\label{t:lb}
For $d \geq 1$ there exists a constant $c_3(d)>0$ such that $\ttop(d,n) \geq c_3(d)(\log n)^{1/d}$ for all $n$.
\end{theorem}

\paragraph{Motivation and background.}
Theorem~\ref{t:pach} is a descendant of the following selection theorem.
  \begin{theorem} [First selection theorem]
\label{t:fs}
Let $P$ be a set of $n$-points in general position in $\R^d$. Then there is a
point in at least $c_4(d) \binom{n}{d+1}$ $d$-simplices spanned by $P$.
\end{theorem}
Theorem~\ref{t:fs} was proved by Boros and F\"{u}redi~\cite{BF84} in the plane
and it was generalized to arbitrary dimension by the first
author~\cite{Barany82}. Relatively recent extensive work of
Gromov~\cite{Gromov10} implies a
topological version of Theorem~\ref{t:fs}; see Theorem~\ref{t:overlap} for
the precise statement of this extension. In addition, Gromov's approach
yielded a significant improvement of the lower bound for the highest possible
value of the constant $c_4(d)$ in Theorem~\ref{t:fs}.

From this point of view, it is desirable to know whether there is a topological
extension of Theorem~\ref{t:pach} which could also possibly be quantitatively
stronger with respect to the constant $c_1(d)$. However,
Theorem~\ref{t:ntp} shows that in the case of this homogeneous selection theorem we
would ask for too much.

\paragraph{A brief proof overview.}
Our proof of Theorem~\ref{t:ntp} partially builds on the approach
from~\cite{Tancer10} where the homogeneous selection theorem was used to
distinguish a geometric and a topological invariant.

For the proof of Theorem~\ref{t:ntp}
we need to exhibit a continuous map
$f\colon (\Delta_{n-1}^{(0)})^{*(d+1)} \to \R^d$ such that $\tau(f)$ is low,
namely at most $c_2(d) n^{1/d}$. Our result is in fact stronger: For some $N \geq (d+1)n$, we construct a map
$f\colon \Delta_{N-1} \rightarrow \R^d$ such that
for {\it any} pairwise disjoint $n$-subsets $V_1,\ldots,V_{d+1}$ of the vertex set of $\Delta_{N-1}$, the restriction of $f$ to $V_1*\cdots*V_{d+1} \cong (\Delta_{n-1}^{(0)})^{*(d+1)}$ satisfies
\begin{equation}
\label{e:strn}
\tau(f_{|V_1*\cdots*V_{d+1}}) \leq c_2(d) n^{1/d}.
\end{equation}

The construction of $f$ proceeds roughly as follows (see Sections
\ref{s:lattice} and \ref{s:maint} for the relevant definitions). Let $L$ be any
finite graded lattice of rank $d+1$ with minimal element $\zht$, whose set of
atoms $A$ satisfies $|A|=N \geq n(d+1)$.  Let $S(A) \cong \Delta_{N-1}$ be
the simplex on the vertex set $A$, and let $\til=L-\{\zht\}$. We first observe
(see Claim \ref{c:nerve}) that there exists a continuous map $g$ from $S(A)$
to the order complex $\Delta(\til)$ such that $g(\langle a_0,\ldots,a_p\rangle)
\subset \Delta(\til_{\leq \vee_{i=0}^p a_i})$ for any atoms $a_0,\ldots,a_p \in
A$ (in words: $\langle a_0,\ldots,a_p\rangle$ maps into the subcomplex below
the join of the atoms $a_0,\ldots,a_p \in A$ in the order complex of $\til$).
Next we define $f:S(A) \rightarrow \R^d$ as the composition $e \circ g$,
where $e: \Delta(\til) \rightarrow \R^d$ is the affine extension of a generic
map from $\til$ to $\R^d$.

Our main technical result, Theorem \ref{t:bfgl},
provides an upper bound on $\tau(f_{|V_1*\cdots*V_{d+1}})$ in terms of the
expansion of the bipartite graph $G_L$ of atoms vs. coatoms of $L$. The desired
bound (\ref{e:strn}) follows from Theorem \ref{t:bfgl} by choosing $L$ to be
the lattice of linear subspaces of the vector space $\FF_q^{d+1}$ over the
finite field with $q$ elements (for suitable $q=q(n,d)$), and utilizing a well
known expansion property of the corresponding graph $G_L$.

The paper is organized as follows: In Section~\ref{s:lattice} we state Theorem~\ref{t:bfgl}
and apply it to prove Theorem~\ref{t:ntp}. The proof of Theorem \ref{t:bfgl} is given in Section
\ref{s:maint}. In Section \ref{s:lbd} we prove Theorem \ref{t:lb} as a direct application of results of Gromov~\cite{Gromov10} and Erd\H{o}s~\cite{Erdos64}.


\paragraph{Subsequent work.}
Considering our work, Bukh and Hubard~\cite{BH17} very recently improved the bound on
$\tau(d,n)$ to $\tau(d,n) \leq 30 (\ln n)^{1/(d-1)}$.

\section{Finite Lattices and Topological Pach Numbers}
\label{s:lattice}

A finite poset $(L,<)$ is a {\it lattice} if for any two element $x,y \in L$ the set
$\{z:z \leq x, z \leq y\}$ has a unique maximal element $x \wedge y$, and the set
$\{z:z \geq x, z \geq y\}$ has a unique minimal element $x \vee y$. In particular, a lattice has a minimal element $\zht$ and a maximal element $\oht$. A lattice $L$ is {\it graded} with rank function $\rk:L \rightarrow \N$, if $\rk(\zht)=0$ and if $\rk(y)=\rk(x)+1$ whenever
$y$ covers $x$ (i.e. $\{z:x \leq z \leq y\}=\{x,y\}$). See Stanley's book \cite{Stanley} for a comprehensive reference on the combinatorics of posets and lattices.
\ \\ \\
Let $L$ be a graded lattice of rank
$\rk(\oht)=d+1$. Let $$A=\{x\in L:\rk(x)=1\}~~,~~C=\{x \in L: \rk(x)=d\}$$
be respectively the sets of {\it atoms} and {\it coatoms} of $L$. For $x \in L$ let $$A_x=\{a \in A : a \leq x\}~~,~~C_x=\{c \in C: x \leq c\}.$$
Let $G_L$ denote the bipartite graph on the vertex set $A \cup C$ with edges
$(a,c) \in A \times C$ iff $a \leq c$. For a set of atoms $Z
\subset A$ let $\Gamma(Z)=\cup_{z \in Z} C_z$ be the neighborhood of $Z$.
\ \\ \\
The main ingredient of the proof of Theorem \ref{t:ntp} is the following connection between $\ttop(d,n)$ and the expansion of $G_L$.
\begin{theorem}
\label{t:bfgl}
Let $L$ be a graded lattice of rank $d+1$ such that $|A| \geq n(d+1)$. Then $m=\ttop(d,n)$ satisfies
\begin{equation*}
\label{e:expl}
\min_{Z \subset A, |Z|=m} |\Gamma(Z)| \leq \frac{d}{d+1} \big(\max_{a \in A} |C_a|+|C|\big).
\end{equation*}
\end{theorem}
\noindent
The proof of Theorem \ref{t:bfgl} is deferred to Section \ref{s:maint}.
\ \\ \\
{\bf Proof of Theorem \ref{t:ntp}:}  Let $n \geq (2d)^d$. By Bertrand's postulate there exists a prime $q$ such that
\begin{equation}
\label{e:bert}
2d \leq \big( (d+1)n \big)^{1/d} \leq q \leq 2\big( (d+1)n \big)^{1/d}.
\end{equation}
Let $\FF_q$ be the finite field of order $q$. Let $L=L(d+1,q)$ denote the
graded lattice of linear subspaces of $\FF_q^{d+1}$ ordered by inclusion, with the natural rank function $\rk(x)=\dim x$ for all $x \in L$. The sets
of atoms and coatoms of $L$ satisfy $|A|=|C|=N_d=\frac{q^{d+1}-1}{q-1}$ and
$|C_a|=N_{d-1}=\frac{q^d-1}{q-1}$ for all $a \in A$. Any two distinct $1$-dimensional subspaces of $\FF_q^{d+1}$ are contained in exactly $N_{d-2}=\frac{q^{d-1}-1}{q-1}$ hyperplanes of $\FF_q^{d+1}$.
Hence, if $a\neq a' \in A$ are two distinct atoms then $$|C_a \cap C_{a'}|=N_{d-2}=\frac{q^{d-1}-1}{q-1}.$$
It follows that if $Z \subset A$, then the family $\{C_a:a \in Z\}$
forms an $N_{d-1}$-uniform hypergraph on vertex set
$\Gamma(Z)$ with $|Z|$ edges, and any two distinct edges intersect in a set of size
$N_{d-2}$. Applying a result of Corr\'{a}di \cite{Corradi69} (see also
exercise 13.13 in \cite{Lovasz93} and Theorem 2.3(ii) in \cite{Alon86}) we obtain
the following lower bound on the expansion of $G_L$.
\begin{equation}
\label{e:expansion}
\begin{split}
|\Gamma(Z)| &\geq \frac{|Z|N_{d-1}^2}{N_{d-1}+(|Z|-1)N_{d-2}}=\frac{|Z|N_{d-1}^2}{q^{d-1}+|Z|N_{d-2}} \\
&=N_d-\frac{q^{d-1}(N_d-|Z|)}{q^{d-1}+|Z|N_{d-2}} \geq
N_d-\frac{q^{d-1}N_d}{|Z|N_{d-2}} \\
&\geq N_d- \frac{q N_d}{|Z|} \geq
N_d- \frac{N_d^{1+\frac{1}{d}}}{|Z|}.
\end{split}
\end{equation}
Next note that (\ref{e:bert}) implies that $|A|=N_d \geq q^d \geq (d+1)n$.
Applying Theorem \ref{t:bfgl} together with (\ref{e:expansion}), it follows that $m=\ttop(d,n)$ satisfies
\begin{equation}
\label{e:fina}
\begin{split}
N_d- \frac{N_d^{1+\frac{1}{d}}}{m} &\leq  \min_{Z \subset A, |Z|=m} |\Gamma(Z)| \\
&\leq \frac{d}{d+1} \big(\max_{a \in A} |C_a|+|C|\big) \\
&=\frac{d}{d+1} (N_{d-1}+N_d).
\end{split}
\end{equation}
The assumption $q \geq 2d$ implies that
\begin{equation}
\label{e:details}
\begin{split}
\frac{N_d}{N_d-dN_{d-1}} &= \frac{q^{d+1}-1}{q^{d+1}-1-d(q^d-1)} \\
&\leq \frac{q^{d+1}}{q^{d+1}-dq^d} = \frac{q}{q-d} \leq 2.
\end{split}
\end{equation}
Rearranging (\ref{e:fina}) and using (\ref{e:details}) and $q^d \leq 2^d (d+1)n$, we obtain
\begin{equation*}
\label{e:ineq1}
\begin{split}
m &\leq \frac{(d+1)N_d^{1+\frac{1}{d}}}{N_d-d N_{d-1}}
\leq 2(d+1)N_d^{\frac{1}{d}} \\
&\leq 2(d+1)\big((d+1)q^d\big)^{1/d} \\
&\leq 2(d+1)\big((d+1)(2^d(d+1)n)\big)^{1/d} \\
&= 4(d+1)\big((d+1)^2 n\big)^{1/d}.
\end{split}
\end{equation*}
{\enp}

\section{Continuous Maps of Finite Lattices}
\label{s:maint}

In this section we prove Theorem \ref{t:bfgl}. We first recall some definitions.
The {\it order complex} $\Delta(P)$ of a finite poset $(P,<)$ is the simplicial complex on the vertex set $P$, whose $k$-simplices are the chains $x_0<\cdots <x_k$ in $P$.

\begin{figure}
\begin{center}
  \scalebox{0.4}{\input{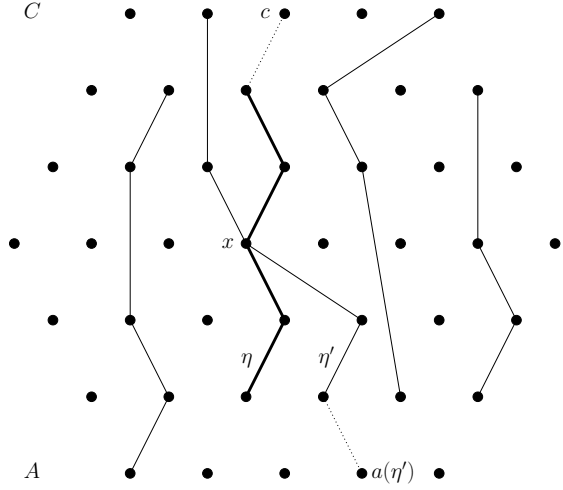}}
  \caption{The bold chain corresponds to $\eta$. The other chains represent simplices of $T'$.}
  \label{figure1}
\end{center}
\end{figure}

Let $L$ be a graded lattice of rank $d+1$ and let $\til=L-\{\zht\}$.
For a subset $\sigma \subset L$ let $\vee \sigma=\vee_{x \in \sigma} x$.
Let $S(A)$ be the simplex on the set $A$ of atoms of $L$ (identified as usual with
its geometric realization).
For $x \in \til$ let $\til_{\leq x}=\{y \in \til: y \leq x\}$.
The main ingredient in the proof of Theorem \ref{t:bfgl} is the following result.
\begin{proposition}
\label{p:main}
There exists a continuous map $f:S(A) \rightarrow \R^d$ such that for any $u \in \R^d$
\begin{equation}
\label{e:mapa}
|\{c \in C: u \in f(\langle A_c \rangle) \}| \leq d \max_{a \in A} |C_a|.
\end{equation}
(Note that, in accordance with our notation, $\langle A_c \rangle$ stands here for the geometric
  realization of $A_c$, considered as a face of $S(A)$.)
\end{proposition}
\noindent
We first note the following
\begin{claim}
\label{c:nerve}
There exists a continuous map $g:S(A) \rightarrow \Delta(\til)$ such that
for all $x \in \til$
\begin{equation*}
\label{c:map}
g(\langle A_x \rangle) \subset \Delta(\til_{\leq x}).
\end{equation*}
\end{claim}
\noindent
{\bf Proof:}
We define $g$ inductively on the $k$-skeleton $S(A)^{(k)}$.  On the vertices
$a \in A$ of $S(A)$ let $g(a)=a$.  Let $0<k \leq |A|-1$ and suppose $g$ has been defined
on $S(A)^{(k-1)}$. Let $\sigma=\langle a_0,\ldots,a_k \rangle \in
S(A)^{(k)}$ and let $y=\vee \sigma$.  For $0 \leq i \leq k$ let
$$\sigma_i=\langle a_0,\ldots,a_{i-1},\widehat{a_i},a_{i+1},\ldots ,a_k \rangle$$
be the $i$-th face of $\sigma$. Let $y_i=\vee \sigma_i$. Then $g$ is defined on
$\sigma_i$ and by induction hypothesis $$g(\sigma_i) \subset \Delta(\til_{\leq
y_i})\subset \Delta(\til_{\leq y}).$$
Being a cone,  $\Delta(\til_{\leq y})$
is contractible and hence $g$ can be continuously extended from the boundary $\partial
\sigma$ to the whole of $\sigma$ so that $g(\sigma) \subset \Delta(\til_{\leq
y})$. It follows in particular that for $x \in \til$
$$
g(\langle A_x \rangle) \subset \Delta(\til_{\leq \vee A_x}) \subset \Delta(\til_{\leq x}).
$$
{\enp}
\noindent {\bf Proof of Proposition \ref{p:main}:} By a general position argument we choose a mapping $e: \til \rightarrow \R^d$ with the following property:
For any pairwise disjoint subsets $S_1,\ldots,S_{d+1} \subset \til$ of cardinalities $|S_i| \leq d$, it holds that $$\bigcap_{i=1}^{d+1} \aff\big(e(S_i)\big) =\emptyset,$$
and thus in particular
\begin{equation}
\label{e:genpos}
\bigcap_{i=1}^{d+1} \relint \conv\big(e(S_i)\big) =\emptyset.
\end{equation}
Extend $e$ by linearity to the whole of $\Delta(\til)$ and let $f=e \circ
g:S(A) \rightarrow \R^d$, where $g$ is the map from Claim \ref{c:nerve}. We claim that the map $f$ satisfies (\ref{e:mapa}).
Let $u \in \R^d$ and let
$$T=\{\eta \in \Delta(\til): u \in \relint e(\langle \eta \rangle)\}.$$
Choose a maximal pairwise disjoint subfamily $T' \subset T$.  It follows by (\ref{e:genpos}) that $|T'| \leq d$.
For each $\eta' \in T'$ choose an atom $a(\eta') \in A$ such that
\begin{equation}
\label{e:atau}
a(\eta') \leq \min \eta'.
\end{equation}
Now let $c \in C$ be such that $u \in f(\langle A_c \rangle)$. Then there exists a
$b \in g(\langle A_c \rangle) \subset \Delta(\til_{\leq c})$ such that
$u=e(b)$. Let $\eta \in T$ be such that
$b \in \relint \langle \eta \rangle$. Then
\begin{equation}
\label{e:tinc}
\eta \in  \Delta(\til_{\leq c}).
\end{equation}
By maximality of $T'$ there exists a simplex $\eta' \in T'$ and a vertex $x \in \eta' \cap \eta$. It follows by (\ref{e:atau}) and (\ref{e:tinc}) that $a(\eta') \leq x \leq c$, i.e. $c \in C_{a(\eta')}$ (see figure \ref{figure1}).
Therefore
$$
|\{c \in C: u \in f(\langle A_c \rangle) \}| \leq \sum_{\eta' \in T'}  |C_{a(\eta')}|
\leq d \max_{a \in A} |C_a|.
$$
{\enp}

\noindent
{\bf Proof of Theorem \ref{t:bfgl}:} Let $L$ be a lattice of rank $d+1$ whose set of atoms $A$ satisfies $|A| \geq (d+1)n$.
Let $V_1,\ldots,V_{d+1}$ be disjoint $n$-subsets of $A$. By Proposition
\ref{p:main} there exists a continuous map $f:S(A) \rightarrow \R^d$ such that for any $u \in \R^d$
\begin{equation*}
\label{e:mapa1}
|\{c \in C: u \in f(\langle A_c \rangle) \}| \leq d \max_{a \in A} |C_a|.
\end{equation*}
Let $m=\ttop(d,n)$. Then there exist $Z_1 \subset V_1,\ldots,Z_{d+1} \subset V_{d+1}$ and a $u \in \R^d$ such that $|Z_i| \geq m$ for all $1 \leq i \leq d+1$ and
$$u \in \bigcap_{z_1\in Z_1,\ldots, z_{d+1} \in Z_{d+1}} f(\langle z_1,\ldots,z_{d+1}\rangle).$$
Write
$$C(Z_1,\ldots,Z_{d+1})=\bigcap_{i=1}^{d+1} \{c \in C: A_c \cap Z_i \neq \emptyset\}.$$
If $c \in C(Z_1,\ldots,Z_{d+1})$ then there exist $z_1 \in Z_1,\ldots,z_{d+1} \in Z_{d+1}$ such that $z_i \leq c$ for all $i$ and hence
$u \in f(\langle z_1,\ldots,z_{d+1}\rangle) \subset f(\langle A_c \rangle)$. Hence by Proposition \ref{p:main}
\begin{equation}
\label{e:uppb}
|C(Z_1,\ldots,Z_{d+1})| \leq d \max_{a \in A} |C_a|.
\end{equation}
On the other hand
\begin{equation}
\label{e:lowb}
\begin{split}
|C(Z_1,\ldots,Z_{d+1})|&=|C-\bigcup_{i=1}^{d+1}(C-\Gamma(Z_i))| \\
&\geq |C|-\sum_{i=1}^{d+1}(|C|-|\Gamma(Z_i)|) =\sum_{i=1}^{d+1}|\Gamma(Z_i)|-d|C| \\
&\geq (d+1)\min_{Z \subset A, |Z|=m} |\Gamma(Z)| -d |C|.
\end{split}
\end{equation}
Theorem \ref{t:bfgl} now follows from (\ref{e:uppb}) and (\ref{e:lowb}).
\enp
\noindent
{\bf Remark:} The mapping $g:S(A) \rightarrow \Delta(\til)$ constructed in Claim \ref{c:nerve} is in general not simplicial. It follows (as of course must be the case by Theorem \ref{t:pach}) that
$f=e \circ g:S(A) \rightarrow \R^d$ is not affine.

\section{The Lower Bound}
\label{s:lbd}

Theorem \ref{t:lb} is a direct consequence of Gromov's topological overlap Theorem \cite{Gromov10} combined with a result of Erd\H{o}s on complete $(d+1)$-partite subhypergraphs in $(d+1)$-uniform dense hypergraphs \cite{Erdos64}. We first recall these results.
Let $X$ be a finite $d$-dimensional pure simplicial complex. For $k \geq 0$, let $f_k(X)=|X(k)|$ denote the number of $k$-dimensional faces of $X$.
Define a positive weight function $w=w_X$ on the simplices of $X$ as follows. For
$\sigma \in X(k)$, let $c(\sigma)=|\{\eta \in X(d): \sigma \subset \eta\}|$
and let $$w(\sigma)=\frac{c(\sigma)}{\binom{d+1}{k+1} f_{d}(X)}.$$
Let $C^k(X)$ denote the space of  $\FF_2$-valued $k$-cochains of $X$ with the coboundary map $d_k:C^k(X) \rightarrow C^{k+1}(X)$.
As usual, the space of $k$-coboundaries is denoted by $d_{k-1}\big(C^{k-1}(X)\big)=B^k(X)$.
For $\phi \in C^k(X)$, let $[\phi]$ denote the image of $\phi$ in
$C^k(X)/B^k(X)$.
Let $$\|\phi\|=\sum_{\sigma \in X(k): \phi(\sigma) \neq 0} w(\sigma)$$
and  $$\|[\phi]\|=\min\{\|\phi+d_{k-1}\psi\|: \psi \in C^{k-1}(X)\}.$$
The {\it $k$-th coboundary expansion constant} of $X$ is
$$
h_k(X)=\min \left\{\frac{\|d_k\phi \|}{\|[\phi]\|}: \phi \in C^k(X)-B^k(X)\right\}.
$$
\noindent
Note that $h_k(X)=0$ iff $\tilde{H}^k(X;\FF_2) \neq 0$. One may regard $h_k(X)$ as a sort of distance between $X$ and the family of complexes $Y$ that satisfy
$\thh^k(Y;\FF_2)\neq 0$. Gromov's celebrated topological overlap result is the following:
\begin{theorem}[Gromov \cite{Gromov10}]
\label{t:overlap}
For any integer $d\geq 0$ and any $\epsilon>0$ there exists a $\delta=\delta(d,\epsilon)>0$ such that if
$h_k(X) \geq \epsilon$ for all $0 \leq k \leq d-1$, then for any continuous map $f:X \rightarrow \R^d$
there exists a point $u \in \R^d$ such that
$$
|\{\sigma \in X(d): u \in f(\sigma)\}| \geq  \delta f_d(X).$$
\end{theorem}
We next describe a result of Erd\H{o}s that generalizes the well known Erd\H{o}s-Stone and K\H{o}v\'ari-S\'{o}s-Tur\'{a}n theorems from graphs to hypergraphs.
\begin{theorem}[Erd\H{o}s \cite{Erdos64}]
\label{t:erdos}
For any $d$ and $c'>0$ there exists a constant $c=c(d,c')>0$ such that for any $(d+1)$-uniform hypergraph $\F$ on $N$-element set $V$ with at least $c' N^{d+1}$ hyperedges,
there exists an $m \geq  c(\log N)^{1/d}$ and disjoint $m$-element sets $Z_1,\ldots,Z_{d+1} \subset V$ such that $\{z_1,\ldots,z_{d+1}\} \in \F$ for all $z_1 \in Z_1,\ldots, z_{d+1} \in Z_{d+1}$.
\end{theorem}
\noindent
{\bf Proof of Theorem \ref{t:lb}:} Recall that $V_1,\ldots,V_{d+1}$ are disjoint $n$-element sets and let
$V=V_1 \cup \cdots \cup V_{d+1}$, $|V|=N=(d+1)n$. Let $X=V_1*\ldots*V_{d+1}$ and let
$f:X \rightarrow \R^d$ be a continuous map. It was shown by Gromov \cite{Gromov10} (see also \cite{DK12,LMM}) that the expansion constants $h_i(X)$ are uniformly bounded away from zero. Concretely, it follows from Theorem 3.3 in \cite{LMM} that $h_i(X) \geq \epsilon=2^{-d}$ for $0 \leq i \leq d-1$. Let $\delta=\delta(d,2^{-d})$. Then by Theorem \ref{t:overlap} there exists a $u\in \R^d$ and a family $\F \subset X(d)$ of cardinality
$$|\F| \geq \delta f_d(X)=\delta n^{d+1}=\delta (d+1)^{-(d+1)} N^{d+1}$$
such that $u \in f(\sigma)$ for all $\sigma \in \F$. Writing $c'= \delta (d+1)^{-(d+1)}$ and $c_3(d)=c(d,c')$, it follows from Theorem \ref{t:erdos} that
there exists an $m \geq c_3(d) (\log N)^{1/d} \geq c_3(d) (\log n)^{1/d}$ and disjoint $m$-sets $Z_1,\cdots,,Z_{d+1} \subset V$ such that $u \in f(\langle z_1,\ldots,z_{d+1}\rangle)$ for all $z_1 \in Z_1,\ldots,z_{d+1} \in Z_{d+1}$. Clearly, there exists a permutation $\pi$ on $\{1,\ldots,d+1\}$ such that $Z_{\pi(i)} \subset V_i$ for all $1 \leq i \leq d+1$.
{\enp}

\ \\ \\
{\bf Acknowledgements}
\\
This research was supported by ERC Advanced Research Grant no 267165 (DISCONV). Imre B\'ar\'any is partially supported by Hungarian National Research Grant K 111827.
Roy Meshulam is partially supported by ISF grant 326/16 and GIF grant 1261/14, Eran Nevo by ISF
grant 1695/15 and Martin Tancer by GA\v{C}R grant 16-01602Y.



\end{document}